\documentclass{amsart}
\usepackage{amssymb}
\usepackage{amscd}
\DeclareFontFamily{OMS}{cmsy}{%
\fontdimen16\font=3pt
\fontdimen17\font=3pt}
\makeatletter
\renewcommand{\subsection}{\@startsection{subsection}{2}{\z@}%
{\baselineskip}{0.5\baselineskip}{\bfseries}}
\makeatother
\newtheorem{thm}{Theorem}[section]
\newtheorem{pro}[thm]{Proposition}
\newtheorem{lem}[thm]{Lemma}
\theoremstyle{definition}

\numberwithin{equation}{section}
\def\dj{d\kern-.30em\raise1.25ex\vbox{\hrule width .3em height .03em}}
\def\Dj{D\rlap{\kern-.70em\raise0.75ex
\vbox{\hrule width .3em height .03em}}}
\newenvironment{pf}{\proof[\proofname]}{\endproof}
\def\Mor{\mathrm{Mor}}
\def\hor{\mathfrak{hor}}
\def\grten{\mathbin{\widehat{\otimes}}}
\def\GP{\cal{C}(P)}
\def\adj{\varpi}
\def\cal{\mathcal}
\def\Bbb{\mathbb}
\def\e{\epsilon}
\def\k{\kappa}
\def\ad{\mathrm{ad}}
\def\id{\mathrm{id}}
\def\gen{\mathrm{gen}}
\def\bla#1{$(${\it #1\/{}}$)$}
\def\gX{\widehat{X}}
\def\GauP{\mathfrak{ad}(P)}
\def\1{\varnothing}
\def\bim#1{\cal{E}_{#1}}
\def\rig{\wp_\Gamma}
\def\inv{{i\!\hspace{0.8pt}n\!\hspace{0.6pt}v}}
\def\im{\mathrm{im}}
\def\Sum{{\displaystyle\sum}}
\def\AgauP{\cal{G}(P)}
\def\WgauP{\widehat{\cal{G}}(P)}
\def\WP{\cal{W}}
\def\BP{\cal{B}}
\def\gphiM{\widehat{\phi}_M}
\def\AdP{\cal{L}}
\def\AdW{\widehat{\cal{L}}}
\def\fWP{\widehat{F}}
\def\fBP{F}
\def\f3{\delta_3}
\def\gf3{\widehat{\delta_3}}
\def\fgau{\Delta}
\def\gfgau{\widehat{\Delta}}
\def\aP{\tau}
\def\gP{\widehat{\tau}}
\def\mM{\mu_M}
\def\sM{\sigma_M}
\def\gsM{\widehat{\sigma}_M}
\def\sAdP{\Sigma}
\def\sAdW{\widehat{\Sigma}}
\def\twst{\chi}
\def\gtwst{\chi}
\def\lt#1{[#1]_1}
\def\rt#1{[#1]_2}
\def\sts{{*}}
\def\filt{\cal{F}}

\begin{document}
\title[Quantum Principal Bundles]{Quantum Gauge Transformations\\
and Braided Structure on\\
Quantum Principal Bundles}

\author{Mi\'co {\Dj}ur{\Dj}evi\'c}
\address{Instituto de Matematicas, UNAM, Area de la Investigacion
Cientifica, Circuito Exterior, Ciudad Universitaria, M\'exico DF, CP
04510, MEXICO}

\begin{abstract}
It is shown that every quantum principal bundle is braided,
in the sense that there exists an intrinsic braid operator twisting
the functions on the bundle. A detailed algebraic analysis of this
operator is performed. In particular, it turns out that the braiding
admits a natural extension to the level of arbitrary differential
forms on the bundle. Applications of the formalism to the study of
quantum gauge transformations are presented. If the structure group
is classical then the braiding is fully compatible with the bundle
structure. This gives a possibility of a direct construction
of the quantum gauge bundle together with its braided quantum group
structure and the corresponding differential calculus. In the general
case, the braiding is not completely compatible with the bundle structure,
however the construction gives a gauge coalgebra containing all the
informations about quantum gauge transformations.
\end{abstract}
\maketitle

\section{Introduction}
\renewcommand{\thepage}{}

In this paper we shall prove that every quantum principal bundle is
equipped with a natural braided structure. This algebraic
structure will be analyzed in details, from various viewpoints.
A particular attention will be given to the study of relations
with quantum gauge transformations, and compatibility with differential
calculus on quantum principal bundles. Our considerations are based on
a general theory of quantum principal bundles, developed in
\cite{d1,d2}.

Quantum principal bundles are noncommutative-geometric \cite{c}
counterparts of classical principal bundles. Quantum groups play
the role of structure groups, and the bundle and the base are
appropriate quantum objects.

In the next section we shall start from a quantum principal bundle $P$
over a quantum space $M$, and introduce a braid operator
$\sM\colon\cal{B}\otimes_M\cal{B}\rightarrow\cal{B}\otimes_M\cal{B}$.
Here $\cal{B}$ is a *-algebra representing $P$ and the tensor product is
over the base space algebra. The map $\sM$ is algebraically
expressed in terms of the right action map $F\colon\cal{B}\rightarrow\cal{B}
\otimes\cal{A}$ and the translation \cite{b-tra} map
$\aP\colon\cal{A}\rightarrow\cal{B}\otimes_M\cal{B}$, where $\cal{A}$ is
a Hopf *-algebra representing the structure group $G$. 
We shall then analyze interrelations between $\sM$ and maps
forming the structure of a quantum principal bundle on $P$.

We shall prove that $\sM$ is fully compatible with the product and
the *-structure on $\cal{B}$. On the other hand, it turns out that $\sM$
is only partially compatible with the action $F$ and the translation $\aP$. 

Section~3 is devoted to the study of relations between the braiding and
differential calculi on the bundle. It turns out that for an arbitrary
differential structure $\Omega(P)$ on the bundle there exists a
natural extension $\gsM\colon\Omega(P)\grten_M\Omega(P)\rightarrow
\Omega(P)\grten_M\Omega(P)$ of the operator $\sM$,
where now the tensor product is over
the graded-differential *-algebra $\Omega(M)$, playing the role of
differential forms on the base in the general theory \cite{d2}. 

The extended operator is constructed with the help of the graded-differential
extension $\gP\colon\Gamma^\wedge\rightarrow\Omega(P)\grten_M\Omega(P)$ of
the translation map, and the pull-back map
$\widehat{F}\colon\Omega(P)\rightarrow\Omega(P)\grten\Gamma^\wedge$. Here
$\Gamma$ is a bicovariant \cite{w3} first-order *-calculus over $G$ and
$\Gamma^\wedge$ is its differential envelope, representing the higher-order
calculus. It is also possible to assume that the complete calculus on
$G$ is based on the associated braided exterior algebra \cite{w3}. 

\renewcommand{\thepage}{\arabic{page}}
We prove that $\gsM$ is completely compatible with the
product and the *-structure on $\Omega(P)$, as well as with the differential
map $d\colon\Omega(P)\rightarrow\Omega(P)$. However $\gsM$ is only partially
compatible with $\widehat{F}$ and $\gP$, because of the incompatibility
appearing already at the level of spaces. 

Then we pass to the appllications of the developed techniques to the
study of quantum gauge transformations. We shall first consider a special
case when the structure group $G$ is classical. The classicality assumption
is equivalent to the disappearance of mentioned incompatibility between
$\sM$ and $\aP,F$. This is further equivalent to the involutivity of $\sM$. 

Quantum principal bundles with classical structure groups give interesting
examples~\cite{d-spin} for the study of quantum phenomenas, but still
incorporating without any change various classical concepts.
In `completely classical' geometry the operator $\sM$ reduces to the
standard transposition.

The full compatibility between $\sM$ and the bundle structure opens a
possibility of constructing in the framework of the formalism
the complete quantum gauge bundle $\GauP$, which is generally understandable
as a braided quantum group \cite{maj} over the quantum space $M$.
If the structure group is classical, this bundle is represented by
the *-algebra $\AdP\subseteq\cal{B}\otimes_M\cal{B}$ consisting of the
elements invariant under the corresponding action of $G$ on
$\cal{B}\otimes_M\cal{B}$.

The compatibility naturally extends to the level of differential structures
if we assume that the calculus on $G$ is classical, too. Then we can conctruct
the appropriate differential calculus on $\GauP$. In classical geometry, the
operator $\gsM$ reduces to the standard graded-transposition. 

In the last section we shall discuss quantum gauge bundles for general
structure groups. In this general context the full quantum gauge
bundle algebra $\AgauP$ can be constructed \cite{d-qgauge} by
applying the structural analysis coming from Tannaka-Krein duality theory
\cite{d-tann}. The present formalism gives a $*$-$\cal{V}$-bimodule
coalgebra $\AdP$, together with a natural fiber-preserving
action map $\fgau\colon\cal{B}\rightarrow
\AdP\otimes_M\cal{B}$.
This is sufficient for the study of the transformation properties.
It turns out
\cite{d-qgauge} that the coalgebra $\AdP$ is naturally
embeddable into the gauge bundle algebra $\AgauP$, via a coalgebra map
$\iota\colon\AdP\rightarrowtail\AgauP$.

Essentially the same picture holds at the graded-differential level. 
Finally, we shall present some explicit calculations with connection forms, 
including gauge transformations of connections, covariant derivative and
curvature. 

\section{Canonical Braiding}

Let us consider a compact \cite{w2} matrix quantum group $G$, represented
by a Hopf *-algebra $\cal{A}$. The elements of $\cal{A}$ play the role of
polynomial functions on $G$. We shall denote by $\phi\colon\cal{A}\rightarrow
\cal{A}\otimes\cal{A}$ the coproduct, and by
$\k\colon\cal{A}\rightarrow\cal{A}$ and $\e\colon\cal{A}\rightarrow\Bbb{C}$
the antipode and the counit map respectively. 

Let $M$ be a quantum space described by a *-algebra $\cal{V}$. Let
$P=(\cal{B},i,F)$ be a quantum principal $G$-bundle over $M$. Here $\cal{B}$
is a *-algebra representing the quantum space $P$, while $F\colon\cal{B}
\rightarrow\cal{B}\otimes\cal{A}$ and $i\colon\cal{V}\rightarrow\cal{B}$ are
*-homomorphisms playing the role of the dualized action of $G$ on $P$ and
the dualized projection of $P$ on $M$ respectively. 

The map $i$ is an isomorphism between $\cal{V}$ and the $F$-fixed-point
subalgebra of $\cal{B}$. In what follows, we shall identify
$\cal{V}$ with its image $\im(i)\subseteq\cal{B}$. 
The algebra $\cal{B}$ is a bimodule over $\cal{V}$,
in a natural manner. We shall use the symbol $\otimes_M$ for the tensor
products over the algebra $\cal{V}$. 

Let $X\colon\cal{B}\otimes_M\cal{B}\rightarrow
\cal{B}\otimes\cal{A}$ be a $\cal{V}$-bimodule homomorphism given by
$$ X(b\otimes q)=bF(q). $$
The freeness property of the action $F$ is equivalent 
to the surjectivity of this map. As explained in \cite{d-ext}, the map $X$ is 
actually bijective. In other words, $P$ is interpretable as
a Hopf-Galois extension.

Let $\aP \colon\cal{A}\rightarrow\cal{B}\otimes_M\cal{B}$ be the
associated translation map \cite{b-tra}. This is a linear map
given by inverting $X$ on $\cal{A}$, that is
$$ \aP (a)=X^{-1}(1\otimes a).$$
Throughout the various computations of this paper,
we shall use a symbolic notation
$$ \aP (a)=\lt{a}\otimes\rt{a}. $$

According to \cite{b-tra} and the analysis of \cite{d2}--Subsection~6.6,
the following identities hold:
\begin{gather*}
\aP (a)^*=\aP [\k(a)^*]\quad\qquad\lt{a}\rt{a}=\e(a)1\\
(\id\otimes F)\aP (a)=\aP(a^{(1)})\otimes a^{(2)}\\
\aP (ac)=\lt{c}\lt{a}\otimes\rt{a}\rt{c}\\
(F\otimes \id)\aP (a)=\lt{a^{(2)}}\otimes\k(a^{(1)})\otimes\rt{a^{(2)}}.
\end{gather*}

Let us also observe that 
\begin{equation}\label{tf=ft}
\aP (a)f=f\aP (a)
\end{equation}
for each $f\in\cal{V}$ and $a\in\cal{A}$. In other words, the elements
$\aP (a)$ are central in the bimodule $\cal{B}\otimes_M\cal{B}$. 

Let $\sM\colon\cal{B}\otimes_M\cal{B}\rightarrow\cal{B}\otimes_M
\cal{B}$ be a linear map defined by
\begin{equation}\label{flip}
\sM(b\otimes q)=\sum_k b_kq\lt{c_k}\otimes \rt{c_k}, 
\end{equation}
where $\Sum_kb_k\otimes c_k=F(b)$. The map $\sM$ is a $\cal{V}$-bimodule 
homomorphism. The left $\cal{V}$-linearity is clear from the definition, while
the right $\cal{V}$-linearity follows from equality \eqref{tf=ft}. Let us
denote by $\mM
\colon\cal{B}\otimes_M\cal{B}\rightarrow\cal{B}$ the product map. 

\begin{pro} \bla{i} The map $\sM$ is bijective, with
\begin{equation}\label{inv}
\sM^{-1}(q\otimes b)=\sum_k\aP \k^{-1}(c_k)qb_k. 
\end{equation}

\smallskip
\bla{ii} The following identity holds
\begin{equation}\label{braid}
(\sM\otimes\id)(\id\otimes\sM)(\sM\otimes\id)
=(\id\otimes\sM)(\sM\otimes\id)(\id\otimes\sM). 
\end{equation}

\smallskip
\bla{iii} The product map $\mM$ is compatible with the braid $\sM$, 
in the sense that
\begin{align}
\sM(\mM\otimes\id)&=(\id\otimes \mM)(\sM\otimes\id)
(\id\otimes\sM)\label{prod-sM1}\\
\sM(\id\otimes \mM)&=(\mM\otimes\id)(\id\otimes\sM)(
\sM\otimes\id). \label{prod-sM2}
\end{align}

\smallskip
\bla{iv} The algebra $\cal{B}$ is $\sM$-commutative. In other words
\begin{equation}\label{comm}
\mM\sM=\mM. 
\end{equation}
\end{pro}
\begin{pf}

A direct computation gives
\begin{multline*}
\sM\sM^{-1}(q\otimes b)=\sum_k\sM\bigl\{\lt{\k^{-1}(c_k)}
\otimes\rt{\k^{-1}(c_k)}qb_k\bigr\}\\=\sum_k
\lt{\k^{-1}(c_k^{(2)})}\rt{\k^{-1}(c_k^{(2)})}qb_k\lt{c_k^{(1)}}\otimes
\rt{c_k^{(1)}}
=\sum_kqb_k\lt{c_k}\otimes\rt{c_k}=q\otimes b, 
\end{multline*}
and similarly
\begin{equation*}
\begin{split}
\sM^{-1}\sM(b\otimes q)&=\sum_k\sM^{-1}\bigl\{b_kq
\lt{c_k}\otimes\rt{c_k}\bigr\}\\
&=\sum_k\lt{\k^{-1}(c_k^{(2)})}\otimes
\rt{\k^{-1}(c_k^{(2)})}b_kq\lt{c_k^{(1)}}\rt{c_k^{(1)}}\\
&=\sum_k \lt{\k^{-1}(c_k)}\otimes\rt{\k^{-1}(c_k)}b_kq=b\otimes q. 
\end{split}
\end{equation*}
We have applied equalities
\begin{align}
\sum_k b_k\lt{c_k}\otimes\rt{c_k}&=1\otimes b\\
\sum_k\aP\k^{-1}(c_k)b_k&=b\otimes 1. 
\end{align}

Let us prove the compatibility relations between $\sM$ and $\mM$. We compute
\begin{equation*}
\begin{split}
\sM(bq\otimes r)&=\sum_{kl}b_kq_lr\lt{d_l}\lt{c_k}\otimes\rt{c_k}\rt{d_l}\\
&=\sum_l(\id\otimes \mM)(\sM\otimes\id)(b\otimes q_lr\lt{d_l}\otimes
\rt{d_l})\\
&=(\id\otimes \mM)(\sM\otimes\id)(\id\otimes \sM)(b\otimes q\otimes r), 
\end{split}
\end{equation*}
where $\Sum_lq_l\otimes d_l=F(q)$. Similarly,
\begin{equation*}
\begin{split}
\sM(b\otimes qr)&=\sum_kb_kqr\lt{c_k}\otimes\rt{c_k}=
\sum_k b_kq\lt{c_k^{(1)}}\rt{c_k^{(1)}}r\lt{c_k^{(2)}}\otimes
\rt{c_k^{(2)}}\\
&=\sum_k(\mM\otimes\id)(\id\otimes\sM)(b_kq\lt{c_k}\otimes\rt{c_k}\otimes r)\\
&=(\mM\otimes\id)(\id\otimes\sM)(\sM\otimes\id)(b\otimes q\otimes r). 
\end{split}
\end{equation*}

Let us now prove that $\sM$ satisfies the braid equation. At first, 
we have
\begin{equation*}
\begin{split}
(\sM\otimes\id)&(\id\otimes\sM)(\sM\otimes\id)(b\otimes q\otimes r)=\\
&=\sum_k(\sM\otimes\id)(\id\otimes\sM)\Bigl(b_kq\lt{c_k}\otimes \rt{c_k}\otimes
r\Bigr)\\&=\sum_k(\sM\otimes\id)\Bigl(b_kq\lt{c_k^{(1)}}\otimes\rt{c_k^{(1)}}
r\lt{c_k^{(2)}}\otimes\rt{c_k^{(2)}}\Bigr)\\
&=\sum_{kl}b_kq_l\lt{c_k^{(3)}}\rt{c_k^{(3)}}r\lt{c_k^{(4)}}\aP \bigl\{
c_k^{(1)}d_l\k(c_k^{(2)})\bigr\}\otimes\rt{c_k^{(4)}}\\
=\sum_{kl}b_kq_lr\lt{c_k^{(3)}}&\lt{\k(c_k^{(2)})}\lt{d_l}\lt{c_k^{(1)}}
\otimes\rt{c_k^{(1)}}\rt{d_l}\rt{\k(c_k^{(2)})}\otimes\rt{c_k^{(3)}}.
\end{split}
\end{equation*}
On the other hand, 
\begin{equation*}
\begin{split}
(\id\otimes\sM)&(\sM\otimes\id)(\id\otimes\sM)(b\otimes q\otimes r)=\\
&=\sum_l(\id\otimes\sM)(\sM\otimes\id)(b\otimes q_lr\lt{d_l}\otimes\rt{d_l})
\\&=\sum_{kl}(\id\otimes\sM)\Bigl\{b_kq_lr\lt{d_l}\lt{c_k}\otimes\rt{c_k}
\otimes\rt{d_l}\Bigr\}\\
&=\sum_{kl}b_kq_lr\lt{d_l}\lt{c_k^{(1)}}\otimes\rt{c_k^{(1)}}\rt{d_l}
\lt{c_k^{(2)}}\otimes\rt{c_k^{(2)}}. 
\end{split}
\end{equation*}
Comparing the obtained expressions, and applying the identity
\begin{equation}
\lt{a^{(2)}}\lt{\k(a^{(1)})}\otimes\rt{\k(a^{(1)})}\otimes\rt{a^{(2)}}
=1\otimes\lt{a}\otimes\rt{a},
\end{equation}
we conclude that \eqref{braid} holds. 

To complete the proof, let us observe that
$$ \mM\sM(b\otimes q)=\sum_k b_kq\lt{c_k}\rt{c_k}=bq=\mM(b\otimes q)$$
in other words, $\cal{B}$ is $\sM$-commutative. 
\end{pf}

In what follows, it will be assumed that $\cal{B}\otimes_M\cal{B}$ is endowed
with the product induced by the operator $\sM$. Explicitly, we define
\begin{equation}\label{prodBB}
(p\otimes b)(q\otimes g)=p\sM(b\otimes q)g. 
\end{equation}
The associativity of this product is ensured by equalities
\eqref{prod-sM1}--\eqref{prod-sM2}. The algebra $\cal{B}\otimes_M\cal{B}$
is unital, with the unit element $1\otimes 1$. The braid equation ensures
that similar $\sM$-induced products can be defined on all $n$-fold
$*$-$\cal{V}$-bimodules
$$ \BP_n=\cal{B}\otimes_M\dots\otimes_M\cal{B}=\cal{B}^{\otimes n}. $$
We shall denote by $\mM^n\colon\BP_n\otimes_M\BP_n\rightarrow\BP_n$
the product map in $\BP_n$. 

Now let us consider the mutual relations between the braid $\sM$ and
the *-structure on $\cal{B}$.

\begin{lem} We have
\begin{equation}
\sM{\sts}={\sts}\sM^{-1}. 
\end{equation}
\end{lem}
\begin{pf} A direct computation gives
\begin{multline*} ({\sts}\sM{\sts})(q\otimes b)=\sum_k{\sts}\bigl(
b_k^*q^*\lt{c_k^*}\otimes\rt{c_k^*}\bigr)=\sum_k\rt{c_k^*}^*\otimes
\lt{c_k^*}^*qb_k\\
=\sum_k \lt{\k^{-1}(c_k)}\otimes \rt{\k^{-1}(c_k)}qb_k=\sM^{-1}(q\otimes b), 
\end{multline*}
where $\Sum_k b_k\otimes c_k=F(b)$. 
\end{pf}

As a direct consequence of the previous lemma, it follows that the formula
\begin{equation}
(b\otimes q)^*=\sM(q^*\otimes b^*)
\end{equation}
defines a new *-involution which is a *-structure on the algebra
$\BP_2$. In a similar way, we
can introduce the *-structure on the higher tensor powers $\BP_n$.
Let us denote by $\twst$ the operators of the standard transposition. 

\begin{lem} \bla{i} The map $\mM\colon\BP_2\rightarrow\cal{B}$
is a *-homomorphism.

\smallskip
\bla{ii} The map $\aP \colon\cal{A}\rightarrow\BP_2$ is a
*-homomorphism, and we have
\begin{equation}
\sM\aP (a)=\aP \k(a)
\end{equation}
for each $a\in\cal{A}$. The following compatibility property holds:
\begin{equation}\label{aP-funct}
(\sM\otimes\id)(\id\otimes\sM)(\aP \otimes\id)=(\id\otimes\aP )\twst.
\end{equation}
\end{lem}

\begin{pf} The hermicity and multiplicativity of $\mM$ directly follow
from the $\sM$-symmetricity property. We have
$$
\sM\aP (a)=\sM(\lt{a}\otimes\rt{a})=\lt{a^{(2)}}\rt{a^{(2)}}
\lt{\k(a^{(1)})}\otimes\rt{\k(a^{(1)})}=\aP \k(a). 
$$
Furthermore,
\begin{multline*}
(\sM\otimes\id)(\id\otimes\sM)(\aP \otimes\id)(a\otimes b)=
(\sM\otimes\id)\Bigl(\lt{a^{(1)}}\otimes\rt{a^{(1)}}b
\lt{a^{(2)}}\otimes\rt{a^{(2)}}\Bigr)\\
=\lt{a^{(2)}}\rt{a^{(2)}}b\lt{a^{(3)}}\lt{\k(a^{(1)})}
\otimes\rt{\k(a^{(1)})}\otimes\rt{a^{(3)}}=b\otimes
\lt{a}\otimes\rt{a}, 
\end{multline*}
which proves \eqref{aP-funct}. Let us check the multiplicativity of
$\aP $. Elementary transformations give
\begin{multline*}
\aP (a)\aP (c)=\lt{a}\sM\bigl(\rt{a}\otimes\lt{c}\bigr)\rt{c}=
\lt{a^{(1)}}\rt{a^{(1)}}\lt{c}\lt{a^{(2)}}\otimes
\rt{a^{(2)}}\rt{c}\\
=\lt{c}\lt{a}\otimes\rt{a}\rt{c}=\aP (ac).
\end{multline*}
Finally, applying the definition of the *-structure in
$\BP_2$ we obtain
$$
\aP(a)^*=\sM\bigl(\lt{\k(a)^*}\otimes \rt{\k(a)^*}\bigr)=
\aP\k[\k(a)^*]=\aP (a^*),
$$
and we conclude that $\aP $ is a *-homomorphism. 
\end{pf}

\begin{lem} We have
\begin{equation}\label{F-funct}
(\sM\otimes\id)(\id\otimes\twst)(F\otimes\id)=(\id\otimes F)\sM.
\end{equation}
\end{lem}

\begin{pf} We compute
\begin{multline*}
(\sM\otimes\id)(\id\otimes\twst)(F\otimes\id)(b\otimes q)=\sum_k
\sM(b_k\otimes q)\otimes c_k=\sum_k b_kq\lt{c_k^{(1)}}\otimes\rt{c_k^{(1)}}
\otimes c_k^{(2)}\\=\sum_k b_k q\lt{c_k}\otimes F\rt{c_k}=
(\id\otimes F)\sM(b\otimes q). \qed
\end{multline*}
\renewcommand{\qed}{}
\end{pf}

For each $n\geq 2$ let us denote by $\fBP_n\colon\BP_n\rightarrow\BP_n\otimes
\cal{A}$ a natural action defined as the $n$-fold direct product of
$F$ with itself. The maps $\fBP_n$ are $*$-$\cal{V}$-bimodule homomorphisms. 

Let $\AdP\subseteq\BP_2$ be a $\cal{V}$-submodule consisting of
$\fBP_2$-invariant elements. This space is closed under the
standard conjugation, however generally it will be not closed under the
$\sM$-induced conjugation. Furthermore, we have a natural projection
map $p_{\AdP}\colon\BP_2\rightarrow\AdP$, given by
\begin{equation}
p_{\AdP}=(\id\otimes h)\fBP_2
\end{equation}
where $h\colon\cal{A}\rightarrow\Bbb{C}$ is the Haar measure \cite{w2}
of $G$. This map is $\cal{V}$-linear and in particular $\AdP$ is a
direct summand in the bimodule $\BP_2$. On the other hand,
$\AdP$ is generally not a subalgebra of $\BP_2$. 

Let us now consider the map $\f3\colon\cal{B}\rightarrow\BP_3$ given by
\begin{equation}\label{f3}
\f3(b)=(\id\otimes\aP)F(b)=\sum_k b_k\otimes\lt{c_k}\otimes\rt{c_k}. 
\end{equation}

\begin{lem} The map $\f3$ is a *-homomorphism and we have
\begin{equation}\label{f3-incl}
\f3(\cal{B})\subseteq \AdP\otimes_M\cal{B}. 
\end{equation}
\end{lem}

\begin{pf} Multiplicativity and hermicity of $\f3$
directly follow from \eqref{aP-funct}. Let us check
the above inclusion. We have
\begin{multline*}
(\fBP_2\otimes\id)\f3(b)=\sum_k\fBP_2(b_k\otimes\lt{c_k})\otimes c_k
=\sum_k b_k\otimes \lt{c_k^{(3)}}\otimes c_k^{(1)}\k(c_k^{(2)})\otimes
\rt{c_k^{(3)}}\\
=\sum_k b_k\otimes \lt{c_k}\otimes 1\otimes\rt{c_k},
\end{multline*}
and we conclude that~\eqref{f3-incl} holds.
\end{pf}

Now let us denote by $\fgau\colon\cal{B}\rightarrow\AdP\otimes_M\cal{B}$ the
map $\f3$ with the restricted codomain. The diagram
\begin{equation}\label{fgau-F}
\begin{CD}
\AdP\otimes_M\cal{B}\otimes\cal{A} @<{\mbox{$\fgau\otimes\id$}}<<
\cal{B}\otimes\cal{A}\\
@A{\mbox{$\id\otimes F$}}AA @AA{\mbox{$F$}}A\\
\AdP\otimes_M\cal{B} @<{\mbox{$\fgau$}}<< \cal{B}
\end{CD}
\end{equation}
is commutative. Indeed, 
\begin{multline*}
(\id^2\otimes F)\f3=(\id^2\otimes F)(\id\otimes \aP)F=
(\id\otimes\aP\otimes\id)(\id\otimes\phi)F\\
=(\id\otimes\aP\otimes\id)(F\otimes\id)F=(\f3\otimes\id)F.
\end{multline*}

\begin{lem}\label{lem:ant} The following identities hold:
\begin{align}
\mM^2(\id^2\otimes\sM)(\f3\otimes\id)&=\mM\otimes 1\\
\mM^2(\sM\otimes\id^2)(\f3\otimes\id)&=1\otimes \mM. 
\end{align}
\end{lem}

\begin{pf}
A direct computation gives
\begin{equation*}
\begin{split} \mM^2(\id^2\otimes\sM)\bigl(\f3(b)\otimes q\bigr)&=
\sum_k\bigl(b_k\otimes\lt{c_k}\bigr)\sM\bigl(\rt{c_k}\otimes q\bigr)\\
&=\sum_kb_k(\id\otimes\mM)(\sM\otimes\id)(\id\otimes\sM)(\aP(c_k)\otimes q)\\
&=\sum_k b_kq\otimes \lt{c_k}\rt{c_k}=bq\otimes 1. 
\end{split}
\end{equation*}
Similarly, we find
\begin{multline*}
\mM^2(\sM\otimes\id^2)\bigl(\f3(b)\otimes q\bigr)=\sum_k
\sM(b_k\otimes \lt{c_k})(
\rt{c_k}\otimes q)\\=\sum_k\bigl\{b_k\lt{c_k^{(2)}}\lt{c_k^{(1)}}
\otimes\rt{c_k^{(1)}}\bigr\}(\rt{c_k^{(2)}}\otimes q)
=\sum_k b_k\lt{c_k}\otimes\rt{c_k}q=1\otimes bq. 
\end{multline*}
We have applied \eqref{aP-funct} and the $\sM$-commutativity of $\cal{B}$. 
\end{pf}

Since the product $\mM$ intertwines $\fBP_2$ and $F$, it follows
that
\begin{equation}
\mM(\AdP)=\cal{V}.
\end{equation}
We shall denote by $\e_M\colon\AdP\rightarrow\cal{V}$ the corresponding
restriction map.
The map $\e_M$ is hermitian, relative to the standard
*-involution on $\BP_2$. 
The intertwining property \eqref{fgau-F} implies
\begin{equation}
(\fgau\otimes\id)(\AdP)\subseteq\AdP\otimes_M\AdP. 
\end{equation}
Let $\phi_M\colon\AdP\rightarrow\AdP\otimes_M\AdP$ be the corresponding
restriction map.

\begin{pro} The following identities hold:
\begin{gather}
(\e_M\otimes\id)\phi_M=(\id\otimes\e_M)\phi_M=\id\label{counit}\\
(\e_M\otimes\id)\fgau=\id\label{e-fgau}\\
(\id\otimes\fgau)\fgau=(\phi_M\otimes\id)\fgau\label{coact}\\
(\phi_M\otimes\id)\phi_M=(\id\otimes\phi_M)\phi_M\label{coasso}.
\end{gather}
\end{pro}

\begin{pf}
We have
\begin{equation*}
(\e_M\otimes\id)\fgau(b)=\sum_k b_k\lt{c_k}\otimes\rt{c_k}=
\sum_k b_k\lt{c_k}\rt{c_k}=b=\sum_k b_k\otimes\lt{c_k}\rt{c_k},
\end{equation*}
which proves \eqref{counit}--\eqref{e-fgau}. Furthermore,
\begin{multline*}
(\phi_M\otimes\id)\fgau(b)=\sum_k b_k\otimes \aP(c_k^{(1)})
\otimes\aP(c_k^{(2)})=\sum_k b_k\otimes\lt{c_k}\otimes
\bigl[F\otimes\aP\bigr]\rt{c_k}\\
=(\id\otimes\fgau)\fgau,
\end{multline*}
which proves \eqref{coact}--\eqref{coasso}. 
\end{pf}

From equalities \eqref{comm} and \eqref{aP-funct} it follows that the map
$X\colon\BP_2\rightarrow\cal{B}\otimes\cal{A}$ is a
*-isomorphism. Actually, the multiplicativity property
of $X$, together with the definition
\eqref{prodBB} of the product in $\BP_2$, completely characterizes the
braiding $\sM$. Moreover, we can define a sequence of *-isomorphisms
$X_n\colon\BP_n\rightarrow\cal{B}\otimes\bigl\{\cal{A}^{\otimes n}\bigr\}$
by equalities
\begin{equation}
X_1=X\qquad X_{n+1}(b\otimes q_n)=(X\otimes\id^n)\bigl(b\otimes X_n(q_n)\bigr).
\end{equation}

We can also introduce generalized translation maps $\aP_n\colon
\cal{A}^{\otimes n}\rightarrow \BP_{n+1}$ by taking the restricted
inverse of $X_n$. It follows that
\begin{equation}
\aP_n(a_1\otimes\dots\otimes a_n)=(\id\otimes\overbrace{\mM\otimes\dots
\otimes\mM}^n\otimes\id)\bigl(\aP(a_1)\otimes\dots\otimes\aP(a_n)\bigr).  
\end{equation}

\section{Extensions to Differential Structures}

In this section we are going to prove that the braid operator $\sM$ is
naturally compatible with an arbitrary differential calculus on a
quantum principal bundle~$P$. All constructions of the previous section will
be incorporated at the level of graded-differential structures.

Let $\Gamma$ be a given bicovariant \cite{w3} first-order *-calculus
over $G$. Let us denote by $\rig$ and $\ell_\Gamma\colon\Gamma
\rightarrow\cal{A}\otimes\Gamma$
be the corresponding right and left action maps. We shall denote by
$\Gamma_{\inv}$ the space of left-invariant elements of
$\Gamma$. The formula
$$\pi(a)=\k(a^{(1)})d(a^{(2)}) $$
defines a projection map $\pi\colon\cal{A}\rightarrow\Gamma_{\inv}$.
Let $\adj\colon\Gamma_{\inv}\rightarrow\Gamma_{\inv}\otimes
\cal{A}$ be the corresponding adjoint action. This map coincides with
the restriction of the right action map on the left-invariant elements. Also,
it is given by projecting the adjoint action of $G$ onto $\Gamma_{\inv}$. 
Explicitly,
$$\adj\pi=(\pi\otimes\id)\ad,\qquad\ad(a)=a^{(2)}\otimes\k(a^{(1)})a^{(3)}.$$

We shall
assume that a higher-order calculus on $G$ is described by the universal
envelope \cite{d1}--Appendix~B of $\Gamma$. The same considerations are
applicable to the braided exterior algebra \cite{w3} associated to $\Gamma$,
relative to the canonical braid operator $\sigma\colon\Gamma_{\inv}^{\otimes 2}
\rightarrow\Gamma_{\inv}^{\otimes 2}$. This map is explicitly given by
$$
\sigma(\eta\otimes\vartheta)=\sum_k\vartheta_k\otimes(\eta\circ c_k), 
$$
where $\Sum_k\vartheta_k\otimes c_k=\adj(\vartheta)$ and $\circ\colon
\Gamma_{\inv}\otimes\cal{A}\rightarrow\Gamma_{\inv}$ is the canonical
right $\cal{A}$-module structure, defined by
$$
\pi(a)\circ b=\pi(ab)-\e(a)\pi(b). 
$$

The Hopf *-algebra structure is naturally liftable from $\cal{A}$ to
a graded differential Hopf *-algebra structure on
$\Gamma^\wedge$. In particular, we have
\begin{gather}
\widehat{\phi}(\vartheta)=\ell_\Gamma(\vartheta)+\rig(\vartheta)\\
\e(\Gamma)=\{0\}\qquad \k[\pi(a)]=-\pi(a^{(2)})\k(a^{(1)})a^{(3)},
\end{gather}
for each $\vartheta\in\Gamma$ and $a\in\cal{A}$. Here $\widehat{\phi}
\colon\Gamma^\wedge\rightarrow\Gamma^\wedge\grten\Gamma^\wedge$ is
the extended coproduct map, and we have denoted by the same symbols
the extended counit and the antipode. 

Let us consider a quantum principal bundle $P=(\cal{B},i,F)$ and let
$\Omega(P)$ be a graded-differential *-algebra representing a differential
calculus on $P$, in the framework of the general theory \cite{d2}.
By definition, $\Omega(P)$ is generated as a differential algebra by
$\cal{B}=\Omega^0(P)$ and there exists a (necessarily unique)
graded-differential homomorphism $\widehat{F}\colon\Omega(P)
\rightarrow\Omega(P)\grten\Gamma^\wedge$ extending the right action map $F$. 
This map is hermitian, and also satisfies
$$
(\id\otimes \widehat{\phi})\widehat{F}=(\widehat{F}\otimes\id)\widehat{F}. 
$$
Let $\hor(P)\subseteq\Omega(P)$ be a graded *-algebra representing
horizontal forms. Horizontal forms are characterized by having trivial
differential properties along vertical fibers, in other words
$$
\hor(P)=\widehat{F}^{-1}\bigl(\Omega(P)\otimes\cal{A}\bigr). 
$$
Finally, a graded-differential *-algebra $\Omega(M)\subseteq\hor(P)$
describing a differential calculus on the base $M$ is defined as a
$\widehat{F}$-fixed point subalgebra of $\Omega(P)$.

There exists a natural right action map $F^\wedge\colon\Omega(P)
\rightarrow\Omega(P)\otimes\cal{A}$, defined as the *-homomorphism
obtained from $\widehat{F}$ by eliminating from the image the
summands with strictly positive second degrees. It turns out \cite{d2}
that $\hor(P)$ is $F^\wedge$-invariant. 

We shall use the symbol $\grten_M$ to denote a tensor product over the
algebra $\Omega(M)$. The map $X$ admits a natural extension $\widehat{X}
\colon\Omega(P)\grten_M\Omega(M)\rightarrow\Omega(P)\grten\Gamma^\wedge$,
given by
$$
\widehat{X}(\varphi\otimes w)=\varphi\widehat{F}(w). 
$$
As explained in \cite{d-ext}, this map is also bijective. This fact allows us
to extend the translation map $\aP\colon\cal{A}\rightarrow\cal{B}\otimes_M
\cal{B}$
to the level of differential algebras, by inverting $\widehat{X}$ on
$\Gamma^\wedge$. In such a way we obtain a map $\gP\colon\Gamma^\wedge
\rightarrow\Omega(P)\grten_M\Omega(P)$. We shall use the same symbolic
notation $\gP(\vartheta)=\lt{\vartheta}\otimes\rt{\vartheta}$ for this
extended map.

Let us define $\Omega(M)$-bimodules $\WP_n$ as $n$-fold tensor products
\begin{equation}
\WP_n=\Omega(P)\grten_M\dots\grten_M\Omega(P), 
\end{equation}
where $n\geq 2$. The differential map $d\colon\Omega(P)\rightarrow
\Omega(P)$ and the *-involution are naturally extendible from $\Omega(P)$
to the spaces $\WP_n$. Furthermore, there exist the natural action maps
$\fWP_n\colon\WP_n\rightarrow\WP_n\grten\Gamma^\wedge$ obtained by taking the
graded-direct products of the action $\widehat{F}$ with itself. The maps
$\fWP_n$ are hermitian, and intertwine the corresponding differentials. 

We shall denote by $\gtwst$ the standard
graded-twist operations. The graded adjoint action map
$\ad\colon\Gamma^\wedge\rightarrow\Gamma^\wedge\grten\Gamma^\wedge$
is given by
$$
\ad(\vartheta)=\gtwst\bigl\{\k(\vartheta^{(1)})\otimes\vartheta^{(2)}\bigr\}
\vartheta^{(3)}. 
$$
The map $\gP\colon\Gamma\rightarrow\WP_2$
possesses the following main algebraic properties:
\begin{gather}
(\id\otimes\widehat{F})\gP(\vartheta)=\gP(\vartheta^{(1)})
\otimes\vartheta^{(2)}\qquad
\fWP_2\gP(\vartheta)=\ad(\vartheta)\\
\gP(\vartheta)^*=\gP(\k(\vartheta)^*)\qquad
\gP[d(\vartheta)]=d\gP(\vartheta)\\
(\widehat{F}\otimes\id)\gP(\vartheta)=\lt{\vartheta^{(2)}}\otimes
\k(\vartheta^{(1)})\otimes\rt{\vartheta^{(2)}}\\
\lt{\vartheta}\rt{\vartheta}=\e(\vartheta)1\qquad
\gP(\vartheta\eta)=(\mM\otimes\id)(\id\otimes\gP)
\gtwst\bigl\{\vartheta\otimes\lt{\eta}\bigr\}\rt{\eta}, 
\end{gather}
where $\mM\colon\WP_2\rightarrow\Omega(P)$ is the product map. 
It turns out that the elements from $\im(\gP)$ graded-commute
with differential forms from $\Omega(M)$.

Let us consider a canonical filtration
$\filt=\Bigl\{\Omega_k(P);k\geq 0\Bigr\}$ of the algebra $\Omega(P)$,
induced by the action $\widehat{F}$ and the grading on $\Gamma^\wedge$.
Let $\gsM\colon\WP_2\rightarrow\WP_2$ be a bimodule homomorphism
defined by
\begin{equation}
\gsM(w\otimes u)=\sum_\alpha
w_\alpha(\mM\otimes\id)(\id\otimes\gP)\gtwst
\bigl\{\vartheta_\alpha\otimes u\bigr\}, 
\end{equation}
where $\Sum_\alpha w_\alpha\otimes\vartheta_\alpha=\widehat{F}(w)$.
By definition, $\gsM$ extends $\sM\colon\BP_2\rightarrow\BP_2$. 
The basic properties of this map are collected in the following

\begin{pro} \bla{i} The map $\gsM$ is bijective, and its inverse is
given by
\begin{equation}\label{g-inv}
\gsM^{-1}(u\otimes w)=\sum_\alpha (\id\otimes\mM)
(\gP\k^{-1}\otimes\id)\gtwst\bigl\{uw_\alpha \otimes\vartheta_\alpha\bigr\}. 
\end{equation}
Furthermore, $\gsM$ is compatible with the filtration $\filt$, in the
sense that 
\begin{equation}\label{gsM-filt}
\gsM\bigl[\Omega_k(P)\grten_M\Omega(P)\bigr]=
\bigl[\Omega(P)\grten_M\Omega_k(P)\bigr],
\end{equation}
for each $k\in\Bbb{N}\cup\{0\}$. In particular
$\gsM\colon\hor(P)\grten_M\Omega(P)\leftrightarrow
\Omega(P)\grten_M\hor(P)$. 

\smallskip
\bla{ii} The map $\gsM$ satisfies the braid equation
\begin{equation}\label{g-braid}
(\gsM\otimes\id)(\id\otimes\gsM)(\gsM\otimes\id)
=(\id\otimes\gsM)(\gsM\otimes\id)(\id\otimes\gsM). 
\end{equation}

\bla{iii} The product map $\mM\colon\WP_2\rightarrow\Omega(P)$ is
compatible with $\gsM$, in a natural manner:
\begin{align}
\gsM(\mM\otimes\id)&=(\id\otimes \mM)(\gsM\otimes\id)
(\id\otimes\gsM)\label{prod-gsM1}\\
\gsM(\id\otimes \mM)&=(\mM\otimes\id)(\id\otimes\gsM)(
\gsM\otimes\id). \label{prod-gsM2}
\end{align}

\bla{iv} The algebra $\Omega(P)$ is $\gsM$-commutative. In other words
\begin{equation}
\mM\gsM=\mM. 
\end{equation}
We also have
\begin{equation}
{*}\gsM=\gsM^{-1}{*}. 
\end{equation}

\bla{v} The braiding $\gsM$ commutes with the differential $d\colon\WP_2
\rightarrow\WP_2$.
\end{pro}

\begin{pf} A direct computation gives
\begin{equation*}
\begin{split}
d\gsM(w\otimes u)&=\sum_\alpha
d(w_\alpha)(\mM\otimes\id)(\id\otimes\gP)\gtwst
(\vartheta_\alpha\otimes u)\\
&\phantom{=}+(-1)^{\partial w}\sum_\alpha
w_\alpha(\mM\otimes\id)(\id\otimes\gP)\gtwst
(\vartheta_\alpha\otimes du)\\
&\phantom{=}+
\sum_\alpha (-1)^{\partial w_{\!\alpha}}
w_\alpha(\mM\otimes\id)(\id\otimes\gP)\gtwst
(d\vartheta_\alpha\otimes u)\\
&=\gsM(dw\otimes u)+(-1)^{\partial w}\gsM(w\otimes
du)=\gsM d(w\otimes u). 
\end{split}
\end{equation*}

Property \eqref{gsM-filt} directly follows from definitions of $\gsM$ and
the filtration $\filt$. The rest of the proof is essentially the same
as for the operator $\sM\colon\BP_2\rightarrow\BP_2$, the only additional
moment is the appearance of graded-twists at the appropriate places. 
\end{pf}

With the help of the operator $\gsM$ we can naturally introduce the *-algebra
structure on the $\Omega(M)$-bimodules $\WP_n$.
Property \bla{v} implies that all differentials
$d\colon\WP_n\rightarrow\WP_n$ are hermitian antiderivations.

Let us consider an $\Omega(M)$-submodule 
$\AdW\subseteq\WP_2$ consisting of all
$\fWP_2$-invariant elements. The space is naturally graded, and we have
$\AdW^0=\AdP$. Furthermore, $\AdW$ is closed under the
standard conjugation and the action of the differential map.
However, $\AdW$ is generally not a subalgebra of $\WP_2$. The product map
induces a bimodule homomorphism $\e_M\colon\AdW\rightarrow\Omega(M)$. This
map satisfies
\begin{align}
\e_M{*}&=*\e_M\\
\e_M d&=d\e_M. 
\end{align}

The map $\fgau\colon\cal{B}\rightarrow\AdP\otimes_M\cal{B}$ admits a natural
extension $\gfgau\colon\Omega(P)\rightarrow\AdW\grten_M\Omega(P)$ given by
\begin{equation}
\gfgau(w)=\sum_\alpha w_\alpha\otimes\lt{\vartheta_\alpha}
\otimes\rt{\vartheta_\alpha}. 
\end{equation}
As a map between $\Omega(P)$ and $\WP_3$, this map is a *-homomorphism. The
map $\gfgau$ preserves the action $\widehat{F}$. In other words, the diagram

\begin{equation}
\begin{CD}
\AdW\grten_M\Omega(P)\grten\Gamma^\wedge @<{\mbox{$\gfgau\otimes\id$}}<<
\Omega(P)\grten\Gamma^\wedge\\
@A{\mbox{$\id\otimes \widehat{F}$}}AA @AA{\mbox{$\widehat{F}$}}A\\
\AdW\grten_M\Omega(P) @<{\mbox{$\gfgau$}}<< \Omega(P)
\end{CD}
\end{equation}
is commutative. The formula
\begin{equation}
\sum_{wu}\gphiM(w\otimes u)=\gfgau(w)\otimes u
\end{equation}
determines a differential bimodule homomorphism $\gphiM\colon\AdW\rightarrow
\AdW\grten_M\AdW$. 

\begin{pro} The following identities hold:
\begin{gather}
(\e_M\otimes\id)\gphiM=(\id\otimes\e_M)\gphiM=\id\\
(\e_M\otimes\id)\gfgau=\id\\
(\id\otimes\gfgau)\gfgau=(\phi_M\otimes\id)\gfgau\\
(\gphiM\otimes\id)\gphiM=(\id\otimes\gphiM)\gphiM.
\end{gather}
In other words, $\AdW$ is a counital differential coalgebra over
$\Omega(M)$ naturally acting on the calculus $\Omega(P)$.
\end{pro}

Let us also mention that all the maps $X_n$ admit natural extensions
$$\gX_n\colon\WP_{n+1}\leftrightarrow
\Omega(P)\grten\underbrace{\Gamma^\wedge\grten\dots\grten\Gamma^\wedge}_n$$
which are isomorphisms of graded-differential *-algebras.
The corresponding partial inverses $\gP_n$ are given by
\begin{equation}
\gP_n(\vartheta_1\otimes\dots\otimes \vartheta_n)=(
\id\otimes\overbrace{\mM\otimes\dots
\otimes\mM}^n\otimes\id)\bigl(\gP(\vartheta_1)\otimes\dots\otimes
\gP(\vartheta_n)\bigr).
\end{equation}

\section{Quantum Gauge Bundles for Classical Structure Groups}
\subsection{The Level of Spaces}

As we have seen in the previous sections, the braid operators $\sM$ and
$\gsM$ are fully compatible with the product maps $\mM$, the *-structure
and the differential $d\colon\Omega(P)\rightarrow\Omega(P)$. Interestingly,
$\sM$ is not completely compatible with the action map
$F\colon\cal{B}\rightarrow\cal{B}\otimes\cal{A}$. It turns out that the
full compatibility between $\sM$ and $F$ holds only in a special case
when $G$ is a classical Lie group. Moreover, differential extensions $\gsM$
and $\widehat{F}$ will be compatible only if the calculus on such a
classical structure group $G$ is assumed to be classical, too.

\begin{pro}
Let $G$ be a compact matrix quantum group, and $P=(\cal{B},i,F)$ a
quantum principal $G$-bundle over a quantum space $M$. Then 
the following conditions are equivalent:

\smallskip
\bla{i} The algebra $\cal{A}$ is commutative. In other words, $G$ is a
classical compact Lie group.

\smallskip
\bla{ii} The following equality holds
\begin{equation}\label{F-wrong}
(\id\otimes\twst)(F\otimes\id)\sM=(\sM\otimes\id)(\id\otimes F). 
\end{equation}

\bla{iii} The following equality holds
\begin{equation}\label{aP-wrong}
(\aP\otimes\id)\gtwst=(\id\otimes\sM)(\sM\otimes\id)(\id\otimes\aP). 
\end{equation}

\bla{iv} The map $\sM$ is involutive.
\end{pro}

\begin{pf} Let us first check equalities \eqref{F-wrong}--\eqref{aP-wrong}.
Direct transformations give
\begin{multline*}
(\id\otimes\sM)(\sM\otimes\id)(b\otimes\aP(a))=
\sum_k(\id\otimes\sM)\bigl[b_k\lt{a}\lt{c_k}\otimes\rt{c_k}
\otimes\rt{a}\bigr]\\=\sum_kb_k\lt{a}\lt{c_k^{(1)}}\otimes
\rt{c_k^{(1)}}\rt{a}\lt{c_k^{(2)}}\otimes\rt{c_k^{(2)}}
\longrightarrow\sum_k b_k\otimes c_k^{(1)}a\otimes c_k^{(2)}, 
\end{multline*}
where $\Sum_k b_k\otimes c_k=F(b)$ and
we have applied the transformation $X_2$ at the end. Similarly,
$$
\aP(a)\otimes b=\lt{a}\otimes\rt{a}\otimes b\longrightarrow \sum_k
b_k\otimes ac_k^{(1)}\otimes c_k^{(2)}. 
$$
Therefore, equality \eqref{aP-wrong} holds if and only if $\cal{A}$ is
commutative. Furthermore,
\begin{multline*}
(\id\otimes\gtwst)(F\otimes\id)\sM(b\otimes q)=
\sum_k(\id\otimes\gtwst)(F\otimes\id)(b_kq\lt{c_k}\otimes\rt{c_k})\\
=\sum_{kl}b_k q_l\lt{c_k^{(3)}}\otimes\rt{c_k^{(3)}}\otimes
c_k^{(1)}d_l\k(c_k^{(2)})\longrightarrow\sum_{kl}
b_kq_l\otimes c_k^{(3)}\otimes c_k^{(1)}d_l\k(c_k^{(2)}),
\end{multline*}
where $\Sum_l q_l\otimes d_l=F(q)$. On the other hand,
$$
(\sM\otimes\id)(\id\otimes F)(b\otimes q)=\sum_{kl} b_kq_l
\lt{c_k}\otimes\rt{c_k}\otimes d_l\longrightarrow\sum_{kl}
b_kq_l\otimes c_k\otimes d_l. 
$$
The obtained expressions coincide if and only if $\cal{A}$ is
commutative. Finally,
let us analyze the question of the involutivity of $\sM$. We compute
\begin{multline*}
\sM^2(b\otimes q)=\sum_k\sM\bigl(b_kq\lt{c_k}\otimes \rt{c_k}\bigr)=
\sum_{kl} b_kq_l\lt{c_k^{(3)}}\rt{c_k^{(3)}}\aP\bigl\{c_k^{(1)}d_l
\k(c_k^{(2)})\bigr\}\\
=\sum_{kl}b_kq_l\aP\bigl\{c_k^{(1)}d_l\k(c_k^{(2)})\bigr\}. 
\end{multline*}
Therefore, $\sM$ will be involutive if and only if
$$ \k(a^{(1)})c a^{(2)}=\e(a)c $$
for each $c,a\in\cal{A}$. This is further equivalent to the commutativity of
$\cal{A}$. 
\end{pf}

It is worth noticing that for involutive braids $\sM$ 
equalities \eqref{F-wrong}--\eqref{aP-wrong} are equivalent to
\eqref{F-funct}--\eqref{aP-funct} respectively. Throughout the rest of this
section, we shall assume that $G$ is a classical compact Lie group. The
bundle $P$ and the space $M$ are arbitrary. 

\begin{pro} The diagram
\begin{equation}\label{tw}
\begin{CD} \BP_2 @>{\mbox{$\fBP_2$}}>> \BP_2\otimes\cal{A}\\
@V{\mbox{$\sM$}}VV @VV{\mbox{$\sM\otimes\id$}}V\\
\BP_2 @>>{\mbox{$\fBP_2$}}> \BP_2\otimes\cal{A}
\end{CD}
\end{equation}
is commutative. In particular, the action $\fBP_2\colon\BP_2\rightarrow
\BP_2\otimes\cal{A}$ is a *-homomorphism and $\AdP$ is a *-subalgebra
of $\BP_2$.
\end{pro}

\begin{pf} We compute
\begin{multline*}
\fBP_2\sM(b\otimes q)=\sum_k\fBP_2\bigl(b_kq\lt{c_k}\otimes\rt{c_k}\bigr)
=\sum_{kl} b_kq_l\lt{c_k^{(3)}}\otimes\rt{c_k^{(3)}}\otimes
c_k^{(1)}d_l\k(c_k^{(2)})c_k^{(4)}\\
=\sum_{kl}b_kq_l\lt{c_k^{(1)}}\otimes\rt{c_k^{(1)}}\otimes c_kd_l
=\sum_{kl}\sM(b_k\otimes d_l)\otimes c_kd_l. \qed
\end{multline*}
\renewcommand{\qed}{}
\end{pf}

The algebra $\AdP$ is covariant relative to the action of the braid
operator $\sM$, in the sense that
\begin{align*}
(\sM\otimes\id)(\id\otimes\sM)\bigl(\AdP\otimes_M\cal{B}\bigr)
&=\cal{B}\otimes_M\AdP\\
(\id\otimes\sM)(\sM\otimes\id)\bigl(\cal{B}\otimes_M\AdP\bigr)&=
\AdP\otimes_M\cal{B}. 
\end{align*}
There exists the induced braiding $\sAdP\colon\AdP\otimes_M\AdP
\rightarrow\AdP\otimes_M\AdP$, given by the restriction of
$$
\sAdP=(\id\otimes\sM\otimes\id)(\sM\otimes\sM)(\id\otimes\sM\otimes\id)\colon
\BP_4\rightarrow\BP_4. 
$$
In particular, we have a natural *-algebra structure in all tensor products
involving $\AdP$ and $\BP_k$. 

It turns out that the braiding $\sAdP$ together with the coalgebra
structure defines a braided quantum group \cite{maj} structure on $\AdP$.
Indeed, the map
$\phi_M\colon\AdP\rightarrow\AdP\grten_M\AdP$ is a *-homomorphism.
Furthermore, 
\begin{equation}
\begin{aligned}\sAdP&=\sAdP^{-1}\\
{*}\sAdP&=\sAdP{*}
\end{aligned}
\qquad
\begin{aligned}
(\id\otimes\phi_M)\sAdP&=(\sAdP\otimes\id)(\id\otimes\sAdP)(\phi_M\otimes\id)\\
(\phi_M\otimes\id)\sAdP&=(\id\otimes\sAdP)(\sAdP\otimes\id)(\id\otimes\phi_M),
\end{aligned}
\end{equation}
as follows from the established full
compatibility between the maps $F$ and $\aP$, and the braiding $\sM$.

From diagram~\eqref{tw} it follows that $\AdP$ is $\sM$-invariant.
Let $\k_M\colon\AdP\rightarrow\AdP$ be the corresponding restriction map. 

\begin{pro} Endowed with maps
$\{\k_M,\e_M,\phi_M, \sAdP\}$ the *-algebra $\AdP$ becomes
a braided-Hopf *-algebra over the quantum space $M$. 
\end{pro}

\begin{pf} 
It is sufficient to check the antipode axiom. However this is a
direct consequence of Lemma~\ref{lem:ant}.
\end{pf}

Geometrically speaking, $\AdP$ represents a quantum gauge bundle $\GauP$
associated to $P$. The elements of $\AdP$ are interpretable as `smooth
functions' on the quantum space $\GauP$. The inclusion $\cal{V}\ni f\mapsto
f\otimes 1\in\AdP$ is interpretable as the dualized fibering of $\GauP$ over
$M$, and $\GauP$ is actually a `bundle of groups'. This is formalized in the
introduced group structure. In particular, the map $\e_M$ is interpretable
as the unit section of the bundle $\GauP$. The homomorphism $\fgau$ is
the dualized left fiberwise bundle action of $\GauP$ on $P$. The diagram
\eqref{fgau-F} expresses the idea that gauge transformations preserve the
structure of the bundle.
In classical geometry $\GauP$ will be the standard
quantum gauge bundle, and actual gauge transformations are naturally
interpretable as smooth sections of $\GauP$.

Following the classical geometry, it is natural to define gauge
transformations as unital $\cal{V}$-linear 
*-homomorphisms $\gamma\colon\AdP\rightarrow\cal{V}$. In our context it
is necessary to impose an additional compatibility condition
\begin{equation}
\gamma(\rho)b=\sum_j b_j\gamma(\rho_j), 
\end{equation}
where
$\Sum_j b_j\otimes \rho_j=(\sM\otimes\id)(\id\otimes\sM)(\rho\otimes b)$.

Let $\GP$ be the set of all maps $\gamma$ of the described type. This set
is a group, in a natural manner. The product and the inverse are given by
\begin{equation}
\gamma\gamma'=(\gamma\otimes\gamma')\phi_M\qquad \gamma^{-1}=\gamma\k^{-1}_M,
\end{equation}
while the unit element is given by the counit map.

The elements of $\GP$ naturally act on $P$, by *-automorphisms of $\cal{B}$
given by the formula
\begin{equation}
\gamma\cdot b=(\gamma\otimes\id)\fgau(b). 
\end{equation}

The following equalities hold
\begin{equation}
(\gamma\gamma')\cdot b=\gamma\cdot(\gamma'\cdot b)\qquad
F(\gamma\cdot b)=\sum_k(\gamma\cdot b_k)\otimes c_k. 
\end{equation}
The group $\GP$ is generally insufficient to cover all symmetry properties
inherent in the action $\fgau$, because $\GauP$ is a quantum object
not reducable to a classical group. 
The presence of the braiding $\sAdP$ reflects the quantum nature of $P$
and $\GauP$.

\subsection{Differential Structures}

The whole reasoning from the previous subsection can be incorporated
at the graded-differential level. Let us assume that $\Gamma$ is the
classical module of differential $1$-forms on $G$. In this case
$\Gamma^\wedge$
gives the standard higher-order calculus on $G$. Let $\Omega(P)$ be an
arbitrary differential calculus on the bundle $P$. We have then

\begin{pro} The following identities hold:
\begin{gather}
(\id\otimes\gtwst)(\widehat{F}\otimes\id)\gsM=(\gsM\otimes\id)(\id\otimes
\widehat{F})\\
\gsM=\gsM^{-1}\\
\fWP_2\gsM=(\gsM\otimes\id)\fWP_2\\
(\gP\otimes\id)\gtwst=(\id\otimes\gsM)(\gsM\otimes\id)(\id\otimes\gP). 
\end{gather}
\end{pro}

We see that $\gP$ and $\widehat{F}$ are completely compatible with the
braiding $\gsM$. In particular, 
the space $\AdW$ is a graded-differential *-subalgebra of $\WP_2$. It is
fully covariant under the action of $\gsM$ and there exists a natural induced
braiding $\sAdW\colon\AdW\grten_M\AdW\rightarrow\AdW\grten_M\AdW$. We have
\begin{equation}
\begin{aligned}
d\sAdW&=\sAdW d\\
\sAdW&=\sAdW^{-1}\\
{*}\sAdW&=\sAdW{*}
\end{aligned}
\qquad
\begin{aligned}
(\id\otimes\gphiM)\sAdW&=(\sAdW\otimes\id)(\id\otimes\sAdW)
(\gphiM\otimes\id)\\
(\gphiM\otimes\id)\sAdW&=(\id\otimes\sAdW)(\sAdW\otimes\id)(\id\otimes\gphiM),
\end{aligned}
\end{equation}

Furthermore, $\gphiM\colon\AdW\rightarrow
\AdW\grten_M\AdW$ is a *-homomorphism. 
The operator $\gsM$
reduces in the space $\AdW$. Let $\k_M\colon\AdW\rightarrow\AdW$ be the
corresponding restriction map.

\begin{pro} Endowed with a system of maps
$\bigl\{\gphiM,\e_M,\k_M,\sAdW\bigr\}$,
the differential *-algebra $\AdW$ becomes a differential braided-Hopf
*-algebra over $M$. \qed
\end{pro}

This gives a natural differential calculus on the quantum gauge bundle. It is
worth mentioning that the differential algebra $\AdW$ is not necessarily
compatible with the classical gauge group $\GP$. This is because
transformations $\gamma\colon\AdP\rightarrow\cal{V}$ from $\GP$ are not
automatically extendible to the appropriate maps on $\AdW$. This
extendibility property can be understood as an additional condition on
the calculus $\Omega(P)$ over the bundle. 

\section{Quantum Gauge Bundles for General Structure Groups}

Let $G$ be an arbitrary compact matrix quantum group, equipped with a
differential calculus $\Gamma$. Let us consider a
quantum principal bundle $P=(\cal{B},i,F)$ over $M$, equipped with
a differential structure $\Omega(P)$. 

The construction of quantum gauge bundles associated to general
quantum principal bundles $P$ can be performed \cite{d-qgauge}
applying the methods developed in \cite{d-tann}. As explained
in \cite{d-tann}, the structure of $P$ is completely encoded in a system
of intertwiner $\cal{V}$-bimodules $\bim{u}=\Mor(u,F)$, associated
to finite-dimensional representations
$u\colon H_u\rightarrow H_u\otimes\cal{A}$ of $G$.
It turns out that the structure of the gauge bundle $\GauP$ is expressible in
terms of the $\cal{V}$-bimodules $\cal{G}_u=\bim{u}\otimes_M\bim{u}^*$
interpretable as consisting of right $\cal{V}$-linear
homomorphisms of $\bim{u}$.

We have the following natural decompositions
$$
\cal{B}=\sideset{}{^\oplus}\sum_{\alpha\in\cal{T}}\cal{B}^\alpha
\qquad \AdP=\sideset{}{^\oplus}\sum_{\alpha\in\cal{T}}\cal{G}_\alpha
$$
of $\cal{V}$-bimodules. Here $\cal{T}$ is a complete set of mutually
inequivalent irreducible representations of $G$. The spaces
$\cal{B}^\alpha\leftrightarrow\bim{\alpha}\otimes H_\alpha$ are the
multiple irreducible bimodules corresponding to the action $F$.

To obtain the full quantum gauge bundle $\GauP$, it is necessary to take into
account all possible bimodules $\cal{G}_u$,
and to factorize through appropriate compatibility relations.
In such a way we obtain a braided quantum group $\AgauP$, together with a
*-coalgebra inclusion $\iota\colon\AdP\rightarrowtail\AgauP$. This map
admits a graded-differential extension $\iota\colon\AdW\rightarrowtail\WgauP$,
where $\WgauP$ is a graded-differential *-algebra describing a complete
calculus \cite{d-qgauge} on $\GauP$. 
In various
interesting (sufficiently regular) special cases the map $\iota$ will be
bijective. This includes,
in particular, locally-trivial quantum principal bundles over classical
smooth manifolds \cite{d1,d-gb}. 

However, for the study of all the phenomenas related to the action
of quantum gauge transformations on the bundle $P$, the coalgebra structure
on $\AdP$ and $\AdW$ is sufficient. We shall analyze in this section
the transformation properties
of basic entities of the general formalism \cite{d2}, under the
action of quantum gauge transformations. 

We have a natural decomposition $\Gamma^{\wedge,\otimes}\leftrightarrow
\cal{A}\otimes\Gamma_{\inv}^{\wedge,\otimes}$, where the left-invariant
part $\Gamma_{\inv}^\wedge\subseteq\Gamma^\wedge$ is the quadratic algebra
explicitly given by
$$ \Gamma_{\inv}^\wedge=\Gamma_{\inv}^\otimes/S_{\inv}^\wedge,
\qquad S_{\inv}^\wedge=\gen(S_{\inv}^{\wedge2})\qquad
S_{\inv}^{\wedge 2}=\Bigl\{\pi(a^{(1)}\otimes\pi(a^{(2)}); a\in\cal{R}\Bigr\},
$$
and $\cal{R}\subseteq\ker(\e)$ is the right $\cal{A}$-ideal
canonically corresponding \cite{w3} to $\Gamma$. In particular, it
follows that the map $\gP\colon\Gamma^\wedge\rightarrow\WP_2$ is
completely determined by its restriction on $\Gamma_{\inv}$.
It turns out that the map
$\gP$ is effectively computed in terms of $\aP$, and the connection forms. 

By definition \cite{d2}, connections on the bundle $P$ are
first-order hermitian linear maps
$\omega\colon\Gamma_{\inv}\rightarrow\Omega(P)$
satisfying the identity
$$
\widehat{F}\omega(\vartheta)=\sum_k\omega(\vartheta_k)\otimes c_k
+1\otimes\vartheta,
$$
where $\Sum_k\vartheta_k\otimes c_k=\adj(\vartheta)$.
If we now fix a grade-preserving splitting
\begin{equation}\label{split-GT}
\Gamma_{\inv}^\otimes\leftrightarrow S_{\inv}^\wedge\oplus
\Gamma_{\inv}^\wedge,
\end{equation}
compatible with  the *-structure and the action $\adj$, then every connection
induces \cite{d2} a natural decomposition
$$
\Omega(P)\leftrightarrow\hor(P)\otimes\Gamma_{\inv}^\wedge\leftrightarrow
\Gamma_{\inv}^\wedge\otimes\hor(P). 
$$
This decomposition plays a fundamental role in various considerations
involving the algebra $\Omega(P)$. Let $\omega$ be an arbitrary
connection on $P$.

\begin{lem} The following identities hold:
\begin{align}
d\aP(a)&=\aP(a^{(1)})\omega\pi(a^{(2)})-\omega\pi(a^{(1)})\aP(a^{(2)})
\label{d-aP}\\
\gP(\vartheta)&=1\otimes\omega(\vartheta)-\sum_k\omega(\vartheta_k)\aP(c_k).
\label{gP-inv}
\end{align}
\end{lem}

\begin{pf} Using the basic transformation property of connections we find
\begin{equation*}
\begin{split}
\gX\Bigl(\aP(a^{(1)})\omega\pi(a^{(2)})-&\omega\pi(a^{(1)})\aP(a^{(2)})\Bigr)
=(1\otimes a^{(1)})\bigl(\omega\pi(a^{(3)})\otimes \k(a^{(2)})a^{(4)}\bigr)\\
&\phantom{=}+1\otimes a^{(1)}\pi(a^{(2)})-\omega\pi(a^{(1)})\otimes a^{(2)}\\
&=1\otimes a^{(1)}\pi(a^{(2)})=1\otimes d(a). 
\end{split}
\end{equation*}
Similarly, we obtain
\begin{multline*}
\gX\Bigl(1\otimes\omega(\vartheta)-\sum_k\omega(\vartheta_k)\aP(c_k)\Bigr)=
1\otimes\vartheta+\sum_k\omega(\vartheta_k)\otimes c_k-
\sum_k\omega(\vartheta_k)X\aP(c_k)\\
=1\otimes\vartheta,
\end{multline*}
which proves \eqref{gP-inv}. 
\end{pf}

Now applying \eqref{gP-inv} and the definition of $\gsM$ we obtain the
following expression for the braiding between connections and arbitrary
differential forms:
\begin{multline}
\gsM\bigl(\omega(\vartheta)\otimes \psi\bigr)=
\sum_k\omega(\vartheta_k)\psi\aP(c_k)-(-1)^{\partial \psi}
\sum_k\psi\omega(\vartheta_k)\aP(c_k)+\\
{}+(-1)^{\partial \psi}\psi\otimes\omega(\vartheta).
\end{multline}

Furthermore, the action of
quantum gauge transformations on connections is given by
\begin{equation}\label{tr-conn}
\gfgau\bigl[\omega(\vartheta)\bigr]=\sum_k\omega(\vartheta_k)\otimes\aP(c_k)
+1\otimes\gP(\vartheta), 
\end{equation}
as directly follows from the definition of the map $\gfgau\colon
\Omega(P)\rightarrow\AdW\grten_M\Omega(P)$, and the basic transformation
property of connections.

Now let us analyze the transformation of the covariant derivative and
the curvature map \cite{d2}. By definition, the curvature of $\omega$ is
a linear map $R_\omega\colon\Gamma_{\inv}\rightarrow\hor(P)$ given by
the structure equation
$$ R_\omega=d\omega-\langle\omega,\omega\rangle $$
where $\langle\ \rangle$ are the brackets associated to
the embedded differential $\delta\colon\Gamma_{\inv}
\rightarrow\Gamma_{\inv}^{\otimes 2}$, coming  from
the splitting~\eqref{split-GT}. More precisely,
$$
\langle\omega,\omega\rangle(\vartheta)=\sum_k\omega(\vartheta_k^1)
\omega(\vartheta_k^2)\qquad \sum_k\vartheta_k^1\otimes\vartheta_k^2
=\delta(\vartheta).
$$

Furthermore, the covariant derivative is a
linear map $D_\omega\colon\hor(P)\rightarrow\hor(P)$ given by
$$ D_\omega(\varphi)=d(\varphi)-(-1)^{\partial\varphi}\sum_k\varphi_k
\omega\pi(c_k) $$
where $F^\wedge(\varphi)=\Sum_k\varphi_k\otimes c_k$.

\begin{lem} The transformation properties of the curvature and covariant
derivative are given by
\begin{align}
\gfgau R_\omega(\vartheta)&
=\sum_k \varsigma(c_k)R_\omega(\vartheta_k)\label{tr-R1}\\
\gfgau D_\omega(\varphi)&
=\sum_k \varsigma(c_k)D_\omega(\varphi_k)\label{tr-D1}, 
\end{align}
where $\varsigma(a)=\lt{\k^{-1}(a^{(1)})}\otimes\aP(a^{(2)})\rt{\k^{-1}
(a^{(1)})}$. 
\end{lem}

\begin{pf} Using the definition of 
of $\gfgau$, and the covariance of the operators $R_\omega$ and $D_\omega$
we obtain
\begin{align}
\gfgau R_\omega(\vartheta)&=\sum_kR_\omega(\vartheta_k)\otimes\aP(c_k)
\label{tr-R2}\\
\gfgau D_\omega(\varphi)&=\sum_kD_\omega(\varphi_k)\otimes\aP(c_k).
\label{tr-D2}
\end{align}
On the other hand,
\begin{multline*}
\sum_kR_\omega(\vartheta_k)\otimes\aP(c_k)
=\sum_k \lt{\k^{-1}(c_k^{(1)})}\otimes\rt{\k^{-1}(c_k^{(1)})}
R_\omega(\vartheta_k)
\aP(c_k^{(2)})\\
=\sum_k \lt{\k^{-1}(c_k^{(1)})}\otimes\aP(c_k^{(2)})\rt{\k^{-1}(c_k^{(1)})}
R_\omega(\vartheta_k). 
\end{multline*}
Similarly, using the covariance of $D_\omega$,
we conclude that equality \eqref{tr-D1} holds. 
\end{pf}

In the definition of the quantum gauge transformation map $\fgau$ and
$\gfgau$ we have not used the braid operators $\sM$ and $\gsM$.
On the other hand, $\gsM$ induces a differential *-algebra structure on
$\WP_3$ such that $\gfgau$ is a differential *-homomorphism. This property
of $\gfgau$ allows us to transform composed algebraic expressions, by
transforming separately their consitutive elements.

As a concrete illustration, let us derive formulas~\eqref{tr-R2}--\eqref{tr-D2}
by a direct computation,
starting from the transformation of connections~\eqref{tr-conn}. We have
$$
\gfgau\langle\omega,\omega\rangle(\vartheta)=
\sum_k\langle\omega,\omega\rangle(\vartheta_k)\otimes\aP(c_k)-
\sum_k\omega(\vartheta_k)\otimes\aP(dc_k)+1\otimes d\gP(\vartheta),
$$
where we have used the identity
$$ \sigma\delta-\delta=(\id\otimes\pi)\adj=c^\top, $$
and the property
\begin{equation}
\langle\gP,\gP\rangle(\vartheta)=\sum_j\gP(\vartheta_j^1)
\gP(\vartheta_j^2)=d\aP(\vartheta).
\end{equation}
Therefore, 
\begin{multline*}
\gfgau R_\omega(\vartheta)=\sum_k d\omega(\vartheta_k)\otimes\aP(c_k)-
\sum_k\omega(\vartheta_k)\otimes \aP(dc_k)+1\otimes d\gP(\vartheta)-\gfgau
\langle\omega,\omega\rangle(\vartheta)\\
=\sum_k d\omega(\vartheta_k)\otimes\aP(c_k)-
\sum_k\langle\omega,\omega\rangle(\vartheta_k)\otimes\aP(c_k)=
\sum_kR_\omega(\vartheta_k)\otimes\aP(c_k). 
\end{multline*}
Furthermore, applying \eqref{d-aP} we obtain
\begin{equation*}
\begin{split}
\gfgau D_\omega(\varphi)&=\gfgau d(\varphi)-(-1)^{\partial\varphi}
\sum_k\bigl(\varphi_k\otimes\aP(c_k^{(1)})\bigr)\gP\pi(c_k^{(2)})\\
&\phantom{=}-(-1)^{\partial\varphi}\sum_k\bigl(\varphi_k\otimes
\aP(c_k^{(1)})\bigr)
\bigl\{\omega\pi(c_k^{(3)})\otimes\aP[\k(c_k^{(2)})c_k^{(4)}]\bigr\}\\
&=\sum_k d(\varphi_k)\otimes\aP(c_k)+(-1)^{\partial\varphi} \varphi_k
\otimes d\aP(c_k)-(-1)^{\partial\varphi}
\sum_k\varphi_k\otimes\gP d(c_k)\\
&\phantom{=}-(-1)^{\partial\varphi}\sum_k\varphi_k\omega\pi(c_k^{(1)})
\otimes\aP(c_k^{(2)})=\sum_k D_\omega(\varphi_k)\otimes \aP(c_k), 
\end{split}
\end{equation*}
and equality \eqref{tr-D2} is proven.

\end{document}